\documentclass[11pt]{article}
\usepackage{amsthm,amsfonts}
\usepackage{amsmath}
\usepackage{amssymb}
\usepackage{graphics,epsfig}
\usepackage{stackrel}
\newtheorem{theorem}{Theorem}[section]

\newtheorem{lemma}{Lemma}[section]

\newtheorem{deft}{Definition}[section]

\newcommand{\argmax}{\operatornamewithlimits{argmax}}
\newcommand{\argmaxproc}{\operatornamewithlimits{argmaxproc}}

\newcommand{\beq}{\begin{equation}}
\newcommand{\enq}{\end{equation}}
\newcommand{\beqn}{$$}
\newcommand{\enqn}{$$}
\newcommand{\beqa}{\begin{eqnarray}}
\newcommand{\enqa}{\end{eqnarray}}
\newcommand{\beqas}{\begin{eqnarray*}}
\newcommand{\enqas}{\end{eqnarray*}}

\newcommand{\bep}{\par\noindent {\sc Proof. }}
\newcommand{\enp}{$\bigtriangleup$\smallskip\par}

\newcommand{\E}{ {\bf E}}
\newcommand{\Prob}{ {\bf P}}

\newtheorem{prop}{Proposition}[section]

\numberwithin{equation}{section}

\begin{document}

\title{Solving optimal stopping problems for L\'{e}vy processes in infinite horizon via $A$-transform }
\author{Elena Boguslavskaya\footnote{elena@boguslavsky.net , Mathematical Sciences, Brunel University, Uxbridge UB8 3PH, UK}\,\,\,\thanks{The author is supported by Daphne Jackson fellowship funded by ESPRC} }
\date{}
\maketitle

\abstract{ We present a method to solve optimal stopping problems in infinite
horizon for a L\'{e}vy process when the reward function can be non-monotone.

To solve the problem we introduce two new objects. Firstly, we  define a random variable $\eta(x)$ which corresponds to the $\argmax$ of the reward function. Secondly, we propose a certain integral transform which can be built on any suitable random variable.  It turns out that this integral transform constructed from $\eta(x)$ and applied to the reward function produces an easy and straightforward description of the optimal stopping rule. We check the consistency of our method with the existing literature, and further illustrate our results with a new example.

The method we propose allows to avoid complicated differential or integro-differential equations which arise if the standard methodology is used. }
\tableofcontents
\section{Introduction}

In recent papers \cite{Surya, Mordecki, MordeckiSalminen, NovikovShiryaevCont, NovikovShiryaevDiscr, KyprianouSurya, MordeckiMishura}  the solutions to optimal stopping problems for L\'{e}vy processes and random walks were found in terms of the maximum/minimum of the process when the solution is one-sided. In \cite{OptStopMarkovSalminen}, the two-sided optimal stopping problem for a strong Markov process was considered.
  In the present paper we extend the results from \cite{Surya} and \cite{NovikovShiryaevCont} to the case of non-monotone reward functions and are able to solve multi-sided optimal stopping problems.\\

In \cite{Surya}, even though some constructions  were defined for a wide class of reward functions, the actual stopping problem was solved only for monotone reward functions.  Paper \cite{OptStopMarkovSalminen} obtains the necessary conditions for a function to solve the two-sided optimal stopping problem. In this paper, we present a constructive method  to solve optimal stopping problems for a fairly general reward function $g$. That is, we show how to find an optimal stopping boundary. The key ingredient here is to construct the integral transform  $\mathcal{A}^{\eta(x)}$ (see Definition \ref{Appellfunction})  on the newly defined random variable $\eta(x) = \argmaxproc_{0 \leq t \leq e_q} g(x+ X_s)  - x$ (see Definition \ref{definition_eta(x)}).
\\

Suppose $X = (X_s)_{s \geq 0}$ is a real-valued L\'{e}vy process. Let $\Prob$ and $\E$
denote the probability and the expectation, respectively, associated with the process
$X$ when started from 0. The natural filtration, generated by $X$, is denoted
by $\mathcal{F} = (\mathcal{F}_t)_{t\geq 0}$, and $\mathcal{M}$ is the set of all stopping times with respect to $\mathcal{F}$.
\\
We aim to find the ``value'' function $V^{*}=V^{*}(\cdot)$ and the optimal stopping time $\tau^*$, such that
\beq
 \label{problemset}V^{*}(x) =
\sup_{\tau \in \mathcal{M}} \E \left(e^{- q \tau}g(x+X_{\tau}) \right) = \E \left(e^{- q \tau^*}g(x+X_{\tau^*}) \right) ,
\enq
where $q>0$, and $g$ is a Borel function to be specified later.
\\
\\
From the general theory of optimal stopping  \cite{ShirOptStop} we know that optimal stopping problems can be solved
through the so-called Markovian method, (see \cite{PeskirShir}, chapter 1). Taking this path, the original problem can be reduced to the corresponding free-boundary problem. Firstly, we have to guess the shape of the free moving boundary (or boundaries). The free moving boundaries divide the  space into subspaces.  We are looking for the optimal stopping boundaries among all possible moving boundaries. The optimal boundaries divide the space into subspaces, namely the ``continuation regions''
 (where it is optimal not to stop, but to continue observations),  and the ``stopping regions'' (where it is optimal to stop the process).
It is not easy and straight forward to guess the shape of the continuation and stopping regions. Additionally, one should always verify that the answer obtained by ``guessing'' is optimal indeed.
\\
\\
 In this paper we show how continuation and stopping regions can be found from the geometrical properties of some function obtained as the integral transform $\mathcal{A}^{\eta(x)}$ of the reward function $g$. We prove that stopping and continuation regions obtained by this method are optimal indeed.
\\
\\
Our algorithm to find the solution to the optimal stopping problem is the following
 \begin{itemize}
 \item We introduce an auxiliary random variable $\eta(x)$ pathwise tracking the value of $X_t$ that achieved the running maximum of $g(x+X)$.
\item We use ${\eta(x)}$ to define the transform $\mathcal{A}^{\eta(x)}$  mapping the reward function $g=g(\cdot)$ into the  function $\mathcal{A}^{\eta(x)}\{g\} (\cdot)$ for each $x$.

\item We define the region $S$ as those arguments at which  $\mathcal{A}^{\eta(x)}\{g\}$(x) is non-negative, i.e.  $S= \{x \mid  \mathcal{A}^{\eta(x)}\{g\}(x)  \geq 0\}$.

\item We define the candidate optimal stopping time as\\ $\tau^+= \inf\{s \geq 0: x+X_s \in S\}$,
 and the candidate value function as \\$V(x) = \E \left( \left[ \mathcal{A}^{\eta(x+ X_{\tau^+})}\{g\}(x+\eta(x)) \right]^+ \right)$.
\item We show that the obtained solution is the optimal solution indeed, i.e the candidate value function $V$ and the candidate optimal stopping time $\tau^+$ coincide with the value function $V^*$ and the optimal stopping time $\tau^*$ from (\ref{problemset}).
\end{itemize}

Finally, we check the consistency of our method with the existing literature by reproducing some well known examples using our method, and then further illustrate our approach by calculating several new examples.

\section{Definitions}
\subsection{${A}$-transform}
Suppose we are given a real function $g=g(\cdot)$ and a random variable
${\nu}$ , with $\E e^{\lambda|\nu|}<\infty$ for some $\lambda>0$. Suppose function $g=g(\cdot)$
has an inverse bilateral Laplace transform $\mathcal{L}^{-1}\{g\}$. For the existence of the inverse bilateral Laplace transform  we can assume, that
$g$ is vanishing at infinity and continuous. Alternatively, we can look at function $g$ as a formal power series, and
take the inverse bilateral Laplace transform formally (for the motivation see Lemma \ref{AppSer}).
\begin{deft} \label{Appellfunction} Let $\nu$ be a random variable, such that $\E e^{\lambda|\nu|}<\infty$ for some $\lambda>0$. The $A$-transform  of function $g=g(\cdot)$ with respect to random variable $\nu$ is a function $Q_g^{\nu}=Q_g^{\nu}(\cdot)$, defined by
\beq    \mathcal{A}^{\nu}\{g\}(y)= Q_g^{\nu}(y) =\int_{-\infty}^\infty \mathcal{L}^{-1}\{g\}(u)  \frac{e^{uy}}{\E e^{u{\nu}}} du . \label{Appell}\enq
\end{deft}
\quad
The  function $Q_g^{\nu}(y) $  is an integral over
the product of the inverse bilateral Laplace transform of function $g$  and
the Esscher transform $\frac{e^{uy}}{\E e^{u \nu}}$ of random variable ${\nu}$.
\\
\\
As it will be shown through examples below, the $A$-transform was designed to convert a reward function $g$
into a function of ``Appell type'', i.e. into a function with properties similar to the Appell function from \cite{NovikovShiryaevCont} and the well-known Appell polynomials. We chose the notation $Q_g^\nu = Q_g^{\nu}(y)$ for the image of   transform $\mathcal{A}^{\nu}$  of  function $g$ in order to be consistent with the existing notation for the Appell function from \cite{NovikovShiryaevCont}, and the Appell polynomials. However, as the term ``Appell function'' is already widely used for an extension of the hypergeometric function to two variables, and the term ``Appell transform'' is used in connection to heat conduction, we decided not to proceed with the term ``Appell'', but to emphasize the ``Appellness'' by denoting the transform by the letter ``${A}$''.\\
\\
One can note, that instead of the bilateral Laplace transform we could have used any exponential transform with the same success. Our choice of the bilateral Laplace is motivated by the desire  to have the Esscher transform in the definition.

\subsection{The random variable $\argmaxproc$ and ${\eta(x)}$}
 By $\varsigma$ we denote the smallest $\argmax_{0 \leq s \leq t}\,\,g (x+X_s)$, in other words  \beq \varsigma=\inf\{s, 0 \leq s \leq t \mid g(x+X_s) \geq g(x+X_u) \mbox{  for any   } u \mbox{   such that   } 0 \leq u \leq t\}.\label{varsigma}\enq
\begin{deft}
\label{definition_eta(x)}
The running $\argmaxproc$ of function $g$, starting at time $0$ and running up to time $t$ is defined by
\beq  \argmaxproc_{0 \leq u \leq t} g(x+X_u) \label{defArgmaxg}=X_{\varsigma},\enq where $\varsigma$ is defined by (\ref{varsigma}).
\end{deft}
We chose the name ``$argmaxproc$'' in order to  emphasize that it is the value of the process on which the $\max$ of the function in question is achieved.
\\
By Definition \ref{definition_eta(x)} we aim to deliver a pathwise construction for the running $\argmaxproc$ of function $g$ over process $X$.
Consider a trajectory of $X$ starting at time $0$, $X_0=0$, and running up to time $t$. Then $g(\argmaxproc_{0 \leq u \leq t} g(x+X_u) )$ is the maximum of the path
\beqn
[0,t]\ni u \rightarrow g(x+X_u) \in \mathbb{R}
\enqn
Note, that if $g$ is a non-decreasing function, then the running $\argmaxproc$ of function $g$ over process $X$ coincides with the running $\max$ of process $X$, i.e.\beq \argmaxproc_{0 \leq u \leq t} g(x+X_u) = \max_{0 \leq u \leq t} (x+X_u).\enq
Similarly, if $g$ is a non-increasing function, then the running $\argmaxproc$ coincides with the running $\min$ of the process
\beq \argmaxproc_{0 \leq u \leq t} g(x+X_u) = \min_{0 \leq u \leq t} (x+X_u).\enq
 Now we are ready to define the random variable $\eta(x)$ which we use in $A$-transform in order to solve our optimal stopping problem.
 \begin{deft}
 Let $e_q$ be an exponentially distributed random variable with mean $1/q$ and independent of the process $X$.
We define the random variable $\eta(x)$  as   \beq \eta(x) = \argmaxproc_{0 \leq s \leq e_{q}}  g(x+X_s) -x \label{defeta}.\enq
\end{deft}
In the same way as above, if $g$ is a non-decreasing function, then $$\eta(x)= \max_{0 \leq s \leq e_{q}}(x+X_s) - x =\max_{0 \leq s \leq e_{q}}X_s,$$
and if $g$ is a non-increasing function, then
$$\eta(x)= \min_{0 \leq s \leq e_{q}}X_s.$$ It is useful to note that  if $g$ is a monotone function, then $\eta(x)$ does not depend on the starting position $x$.

\section{Main results. Solution to the optimal stopping problem}
\subsection{The candidate value function and the candidate optimal stopping time}
To solve our optimal stopping problem we have to find the value function $V^{*}=V^{*}(\cdot)$ and the optimal stopping time $\tau^*$, such that
\beq
 \label{problemset}V^{*}(x) =
\sup_{\tau \in \mathcal{M}} \E \left(e^{- q \tau}g(x+X_{\tau}) \right) = \E \left(e^{- q \tau^*}g(x+X_{\tau^*}) \right) ,
\enq
Let us introduce {\it a candidate optimal stopping time} $\tau^+$ and a {\it candidate value function} $V$.
\\
Define  set $S$ as the set of all $x$ such that $A$-transform of $g$ with respect to $\eta(x)$ is non-negative
\beq \label{definitionS} S:=\{x \mid \mathcal{A}^{\eta(x)}\{g\}(x) \geq 0\},\enq
where $\eta(x)$ is defined by (\ref{defeta}).
\\
Let the candidate optimal stopping time $\tau^+$ be defined as the first moment at which $x+X_s$ reaches $S$
\beq \tau^+:= \inf\{s \geq 0: x+ X_s \in S\}.\label{definition_tau+}\enq
Now define the candidate value function $V$ as
\beq
\label{definition_V}
V(x) := \E \left( \left[\mathcal{A}^{\eta(x+X_{\tau^+})}\{g\}(x+\eta(x)) \right]^+ \right).
\enq

\subsection{Auxiliary lemmas}
To prove the optimality of the candidate solution we need the following lemmas.
\begin{lemma} Suppose $g=g(y)$ and the image of $A$-transform of $g$ with respect to $\eta(x)$, i.e. $\mathcal{A}^{\eta(x)}\{g\}(y)$,  are co-monotone functions in $y$ for each fixed $x$ and for those $y$ where $\mathcal{A}^{\eta(x)}\{g\}(y) \geq 0$. Then for any stopping moment $\tau$ and for any $x$ we have $$\E \left( \left[\mathcal{A}^{\eta(x+X_{\tau^+})}\{g\}(x+\eta(x)) \right]^+ \right) \geq  \E \left(e^{- q \tau}g(x+X_{\tau}) \right).$$ \label{lemma ineq}
\end{lemma}
\bep
Let $\tilde{X}$ be an independent version of process $X$, and let
$$
\tilde{X}_\sigma = \tilde{\eta}(x + X_{\tau^+})=\argmaxproc_{0 \leq s \leq e_q} g\left(x + X_{\tau^+} + \tilde{X}_s\right) -x - X_{\tau^+}.
$$
Now, using the tower property of conditional expectation, definitions of $\eta(x)$ and $\tilde{\eta}(x)$, and the co-monotonicity of function $g$ and its $A$-transform, we have the following chain of equalities and inequalities
for any stopping moment $\tau$ and any $x$
\beqas
\E \left( \left[\mathcal{A}^{\eta(x+X_{\tau^+})}\{g\}(x+ \eta(x))\right]^+\right) &=&
\E \left(\E\left(\left[\mathcal{A}^{\eta(x+X_{\tau^+})}\{g\}(x+ \eta(x))\right]^+\mid \mathcal{F}_{\tau}\right)\right)\\
&=&
\E \left(\E\left(\left[\mathcal{A}^{\eta(x+X_{\tau^+})}\{g\}(\argmaxproc_{0 \leq t \leq e_q} g\left(x + X_t\right))\right]^+\mid \mathcal{F}_{\tau}\right)\right)\\
&=& \E\left(\E\left(\sup_{0 \leq t \leq e_q}\left(\left[\mathcal{A}^{\eta(x+X_{\tau^+})}\{g\}(x+ X_t)\right]^+\right)\mid \mathcal{F}_{\tau}\right) \right)\\
&\geq& \E\left(\E\left(\sup_{\tau \leq t \leq e_q}\left(\left[\mathcal{A}^{\eta(x+X_{\tau^+})}\{g\}(x+ X_t)\right]^+\right)
\mid \mathcal{F}_{\tau}\right)\mathbf{1}_{\{e_q > \tau\}}\right)\\
&=& \E\left(\E\left(\sup_{\tau \leq t \leq e_q}\left(\left[\mathcal{A}^{\eta(x+X_{\tau^+})}\{g\}(x+X_\tau+ X_t -X_\tau)\right]^+\right)
\mid \mathcal{F}_{\tau}\right)\mathbf{1}_{\{e_q > \tau\}}\right)\\
%
&=& \E\left(\E\left(\sup_{0\leq s \leq e_q}\left(\left[\mathcal{A}^{\eta(x+X_{\tau^+})}\{g\}(x+X_\tau+ \tilde{X_s})\right]^+
\right)
\mid \mathcal{F}_{\tau}
\right)\mathbf{1}_{\{e_q > \tau\}}\right)\\
&\geq&
\E\left(\E\left(\left(\left[\mathcal{A}^{\eta(x+X_{\tau^+})}\{g\}(x+X_\tau+ \tilde{X}_\sigma)\right]^+
\right)
\mid \mathcal{F}_{\tau}
\right)\mathbf{1}_{\{e_q > \tau\}}\right)\\
&=& \E\left(\E\left(\left(\left[\mathcal{A}^{\eta(x+X_{\tau^+})}\{g\}(x+X_\tau+ \tilde{\eta}(x+X_{\tau^+}))\right]^+
\right)
\mid \mathcal{F}_{\tau}
\right)\mathbf{1}_{\{e_q > \tau\}}\right)\\
&\geq& \E\left(\E\left(\left(\mathcal{A}^{\eta(x+X_{\tau^+})}\{g\}(x+X_\tau+ \tilde{\eta}(x+X_{\tau^+}))
\right)
\mid \mathcal{F}_{\tau}
\right)\mathbf{1}_{\{e_q > \tau\}}\right)\\
&=&
\E\left(\E\left(\int_{-\infty}^\infty \mathcal{L}^{-1}(g)(u)\frac{e^{u(x+X_{\tau}+\tilde{\eta}(x+ X_{\tau^+}))}}{\E e^{u\eta(x+ X_{\tau^+})}}du
\mid \mathcal{F}_{\tau}
\right) \mathbf{1}_{\{e_q > \tau\}}\right)\\
&=& \E\left(\int_{-\infty}^\infty \mathcal{L}^{-1}(g)(u)\frac{e^{u(x+ X_{\tau})} \E (e^{u\tilde{\eta}(x+ X_{\tau^+})}
\mid \mathcal{F}_{\tau}
)}{\E e^{u\eta(x+ X_{\tau^+})}}du \mathbf{1}_{\{e_q > \tau\}}\right)\\
&=& \E\left(\int_{-\infty}^\infty \mathcal{L}^{-1}(g)(u)\frac{e^{u(x+ X_{\tau})} \E (e^{u\tilde{\eta}(x+ X_{\tau^+})})}{\E e^{u\eta(x+ X_{\tau^+})}}du \, \mathbf{1}_{\{e_q > \tau\}}\right)\\
&=& \E \left(\int_{-\infty}^\infty \mathcal{L}^{-1}(g)(u)e^{u(x+ X_{\tau})} du\, \mathbf{1}_{\{e_q > \tau\}}\right)\\
&=& \E \left(g(x+X_{\tau} )e^{-q \tau}\right).
\enqas
\enp
\begin{lemma} Suppose $g=g(y)$ and the image of $A$-transform of $g$ with respect to $\eta(x)$, i.e. $\mathcal{A}^{\eta(x)}\{g\}(y)$,  are co-monotone functions in $y$ for each fixed $x$ and for those $y$ where $\mathcal{A}^{\eta(x)}\{g\}(y) \geq 0$. Let  $\tau^+$ be defined by (\ref{definition_tau+}). Then for  any $x$ we have $$\E \left( \left[\mathcal{A}^{\eta(x+X_{\tau^+})}\{g\}(x+\eta(x)) \right]^+ \right) =  \E \left(e^{- q {\tau^+}}g(x+X_{\tau^+}) \right).$$ \label{lemma eq}
\end{lemma}
\bep
Let $\tilde{X}$ be an independent version of process $X$, and let
$$
\tilde{\eta}(x + X_{\tau^+})=\argmaxproc_{0 \leq s \leq e_q} g\left(x + X_{\tau^+} + \tilde{X}_s\right) -x - X_{\tau^+}.
$$
Now, using the tower property of conditional expectation, definitions of $\eta(x)$ and $\tilde{\eta}(x)$, and the co-monotonicity of function $g$ and its $A$-transform we have the following chain of equalities for any $x$
\beqas
\E \left( \left[\mathcal{A}^{\eta(x+X_{\tau^+})}\{g\}(x+ \eta(x))\right]^+\right) &=&
\E \left(\E\left(\left[\mathcal{A}^{\eta(x+X_{\tau^+})}\{g\}(x+ \eta(x))\right]^+\mid \mathcal{F}_{{\tau^+}}\right)\right)\\
&=&
\E \left(\E\left(\left[\mathcal{A}^{\eta(x+X_{\tau^+})}\{g\}(\argmaxproc_{0\leq t \leq e_q} g\left(x + X_t\right))\right]^+\mid \mathcal{F}_{{\tau^+}}\right)\right)\\
&=& \E\left(\E\left(\sup_{0 \leq t \leq e_q}\left(\left[\mathcal{A}^{\eta(x+X_{\tau^+})}\{g\}(x+ X_t)\right]^+\right)\mid \mathcal{F}_{{\tau^+}}\right) \right)\\
&=& \E\left(\E\left(\sup_{{\tau^+} \leq t \leq e_q}\left(\left[\mathcal{A}^{\eta(x+X_{\tau^+})}\{g\}(x+ X_t)\right]^+\right)
\mid \mathcal{F}_{{\tau^+}}\right)\mathbf{1}_{\{e_q > {\tau^+}\}}\right)\\
&=& \E\left(\E\left(\sup_{{\tau^+} \leq t \leq e_q}\left(\left[\mathcal{A}^{\eta(x+X_{\tau^+})}\{g\}(x+X_{\tau^+}+ X_t -X_{\tau^+})\right]^+\right)
\mid \mathcal{F}_{{\tau^+}}\right)\mathbf{1}_{\{e_q > {\tau^+}\}}\right)\\
%
&=& \E\left(\E\left(\sup_{0\leq s \leq e_q}\left(\left[\mathcal{A}^{\eta(x+X_{\tau^+})}\{g\}(x+X_{\tau^+}+ \tilde{X_s})\right]^+
\right)
\mid \mathcal{F}_{{\tau^+}}
\right)\mathbf{1}_{\{e_q > {\tau^+}\}}\right)\\
&=&\E\left(\E\left(\left[\mathcal{A}^{\eta(x+X_{\tau^+})}\{g\}\left(\argmaxproc_{0 \leq s \leq e_q} g\left(x+X_{\tau^+}+ \tilde{X_s}\right)\right)\right]^+\,
\mid \mathcal{F}_{{\tau^+}}
\right)\mathbf{1}_{\{e_q > {\tau^+}\}}\right)\\
&=& \E\left(\E\left(\left[\mathcal{A}^{\eta(x+X_{\tau^+})}\{g\}(x+X_{\tau^+}+ \tilde{\eta}(x+X_{{\tau^+}}))\right]^+
\mid \mathcal{F}_{{\tau^+}}
\right)\mathbf{1}_{\{e_q > {\tau^+}\}}\right)\\
&=& \E\left(\E\left(\mathcal{A}^{\eta(x+X_{\tau^+})}\{g\}(x+X_{\tau^+}+ \tilde{\eta}(x+X_{\tau^+}))
\mid \mathcal{F}_{{\tau^+}}
\right)\mathbf{1}_{\{e_q > {\tau^+}\}}\right)\\
&=&
\E\left(\E\left(\int_{-\infty}^\infty \mathcal{L}^{-1}(g)(u)\frac{e^{u(x+X_{{\tau^+}}+\tilde{\eta}(x+ X_{{\tau^+}}))}}{\E e^{u\eta(x+ X_{\tau^+})}}du
\mid \mathcal{F}_{{\tau^+}}
\right) \mathbf{1}_{\{e_q > {\tau^+}\}}\right)\\
&=& \E\left(\int_{-\infty}^\infty \mathcal{L}^{-1}(g)(u)\frac{e^{u(x+ X_{{\tau^+}})} \E (e^{u\tilde{\eta}(x+ X_{\tau^+})}
\mid \mathcal{F}_{{\tau^+}}
)}{\E e^{u\eta(x+ X_{\tau^+})}}du \mathbf{1}_{\{e_q > {\tau^+}\}}\right)\\
&=& \E\left(\int_{-\infty}^\infty \mathcal{L}^{-1}(g)(u)\frac{e^{u(x+ X_{{\tau^+}})} \E (e^{u\tilde{\eta}(x+ X_{\tau^+})})}{\E e^{u\eta(x+ X_{\tau^+})}}du \, \mathbf{1}_{\{e_q > {\tau^+}\}}\right)\\
&=& \E \left(\int_{-\infty}^\infty \mathcal{L}^{-1}(g)(u)e^{u(x+ X_{{\tau^+}})} du\, \mathbf{1}_{\{e_q > {\tau^+}\}}\right)\\
&=& \E \left(g(x+X_{{\tau^+}} )e^{-q {\tau^+}}\right).
\enqas
\enp

\subsection{The main theorem}
\begin{theorem} Suppose $g=g(y)$ and the image of $A$-transform of $g$ with respect to $\eta(x)$, i.e. $\mathcal{A}^{\eta(x)}\{g\}(y)$,  are co-monotone functions in $y$ for each fixed $x$ and for those $y$ where $\mathcal{A}^{\eta(x)}\{g\}(y) \geq 0$. Let $\tau^+$ and $V=V(x)$ be defined respectively by (\ref{definition_tau+}) and (\ref{definition_V}).  Then the stopping time $\tau^+$  and the function $V$ are the optimal stopping time and the value function for the problem (\ref{problemset}), i.e. $\tau^*=\tau^+$ and  $V^*(x)=V(x)$.
\end{theorem}
\bep
Indeed, by lemma \ref{lemma ineq}  $V(x) \geq \E \left( e^{-q\tau} g(x+X_{\tau})\right) $ for {\it any} stopping time $\tau$.
Naturally, it also holds for the optimal stopping time $\tau^*$ $$V(x) \geq \E \left( e^{-q\tau^*} g(x+X_{\tau^*})\right)=\sup_\tau \E \left( e^{-q\tau} g(x+X_{\tau})\right)  = V^{*}(x).$$
On the other hand,
 for $\tau^+$ we have
 \beqas
V^{*}(x) &=& \sup_{\tau} \E \left( e^{-q\tau} g(x+X_{\tau}) \right) \\ &\geq& \E \left( e^{-q \tau^+} g(x+X_{\tau^+}) \right)
 \stackrel{lemma \ref{lemma eq}}{=} V(x).
\enqas
Thus,  $V(x) \geq V^{*}(x) \geq V(x)$ for any $x$, i.e. $V(x) = V^{*}(x)$ for any $x$.
Furthermore, $V(x)= \E \left(e^{-q \tau^+} g( x+X_{\tau^+}) \right)$ by lemma \ref{lemma eq}. Thus, $\tau^+$ is the optimal stopping time.
\enp

\section{Properties of $A$-transform}
\subsection{The averaging property}

\begin{lemma} Let $\nu$ be a random variable with $\E e^{\lambda|\nu|}<\infty$ for some $\lambda>0$, and $\mathcal{A}^\nu\{g\}(y)$ be an ${A}$-transform  of function $g$ with respect to random variable $\nu$ given by Definition \ref{Appellfunction}. Then transform $\mathcal{A}^\nu$ satisfies the averaging property \label{averagingproperty}
$\E(\mathcal{A}^\nu\{g\}(y+\nu))= g (y) \label{EQ}. $
\end{lemma}
\bep
Indeed,
\beqn \E(\mathcal{A}^\nu\{g\}(y+\nu))= \int_{-\infty}^\infty \mathcal{L}^{-1}\{g\}(u)
\frac{\E e^{u{(y+\nu)}}}{\E e^{ u {\nu}}} du = \int_{-\infty}^\infty \mathcal{L}^{-1}\{g\}(u)
e^{uy} du = g(y). \enqn \enp

\subsection{The martingale property of $\left(\mathcal{A}^{X_t}\{g\}(X_t)\right)_{t \geq 0}$ for a L\'{e}vy process $X$. }
\begin{lemma} Let $\left(X_t\right)_{t\geq0}$ be a L\'{e}vy process such that there is some $\lambda>0$ such that $\E e^{\lambda|X_t|}<\infty$, $t \geq0$. Then $\left(\mathcal{A}^{X_t}\{g\}(X_t)\right)_{t\geq0}$ is a martingale.
\end{lemma}
\bep Indeed, for all $t$ we have $\E|\mathcal{A}^{X_t}\{g\}(X_t) | < \infty$ and
\beqas
\E \left( \left.\mathcal{A}^{X_t}\{g\}(X_t) \right|\mathcal{F}_s\right)
&=&\E\left( \left.\int_{-\infty}^\infty \mathcal{L}^{-1}\{g\}(u)\frac{ e^{u{X_t}}}{\E e^{ u X_t}}du \right| \mathcal{F}_s \right)
\\
&=&\int_{-\infty}^\infty \mathcal{L}^{-1}\{g\}(u)\E\left( \left.\frac{ e^{u{X_t}}}{\E e^{ u X_t}} \right| \mathcal{F}_s \right)du\\
&=& \int_{-\infty}^\infty \mathcal{L}^{-1}\{g\}(u) \frac{ e^{u{X_s}}}{\E e^{ u X_s}} du\\
&=& \mathcal{A}^{X_s}\{g\}(X_s).
\enqas
\enp
\subsection{The linearity}
The linearity of ${A}$-transform  follows from linearity of the inverse bilateral Laplace transform.
\begin{lemma}Suppose ${A}$-transform exists for real functions $f$ and $g$.
Then
\beq
\mathcal{A}^\nu \{c_1f+ c_2g\}(y) = c_1\mathcal{A}^\nu \{f\}(y)+c_2\mathcal{A}^\nu\{g\}(y),
\enq
where $c_1$ and $c_2$ are some constants.
\end{lemma}
\bep
\beqas
\mathcal{A}^\nu\{c_1f+ c_2g\}(y) &=& \int_{-\infty}^{\infty}\mathcal{L}^{-1}\{c_1f+c_2g\}(u)\frac{e^{uy}}{\E e^{u\nu}}du\\
&=& \int_{-\infty}^{\infty}\left(c_1\mathcal{L}^{-1}\{f\}(u)+c_2\mathcal{L}^{-1}\{g\}(u)\right)\frac{e^{uy}}{\E e^{u\nu}}du\\
&=& c_1\int_{-\infty}^{\infty}\mathcal{L}^{-1}\{f\}(u)\frac{e^{uy}}{\E e^{u\nu}}du + c_2\int_{-\infty}^{\infty}\mathcal{L}^{-1}\{g\}(u)\frac{e^{uy}}{\E e^{u\nu}}du\\
&=&c_1\mathcal{A}^\nu\{f\}(y)+c_2\mathcal{A}^\nu\{g\}(y).
\enqas
\enp
\subsection{The differential property}
\begin{lemma} Let $f$ be a differentiable function on $\mathbb{R}$. Suppose ${A}$-transform with respect to $\nu$ of $f$    is differentiable. Then it satisfies
the following differential property \beq
\frac{d}{dy} \left(\mathcal{A}^\nu\{f\}(y)\right)=\mathcal{A}^\nu\left\{ \frac{df}{dx} \right\}(y).  \label{DifferentialPropertyOfAppellFunction}
 \enq
\end{lemma}
\bep It is well known  that  if $\mathcal{L}^{-1}(f)$ is the inverse
bilateral Laplace transform of $f(y)$, then $u\mathcal{L}^{-1}(f)$ is the inverse bilateral Laplace
transform of $\frac{d}{dy}f(y)$. Therefore, we have

\beqas \frac{d}{dy} \left(\mathcal{A}^\nu\{f\}(y)\right)= \frac{d}{dy} Q_f^\nu(y) &=& \int_{-\infty}^\infty \mathcal{L}^{-1}(f)(u)
\frac{d}{dy}\left(\frac{e^{uy}}{\E e^{ u \nu}}\right) du
\\&=& \int_{-\infty}^\infty \left(u \mathcal{L}^{-1}(f)(u)\right)
\frac{e^{uy}}{\E e^{ u \nu}} du \\
&=& \int_{-\infty}^\infty \mathcal{L}^{-1}\left(\frac{df}{dx}\right)(u)
\frac{e^{uy}}{\E e^{ u \nu}} du
\\&=&  \mathcal{A}^\nu\left\{ \frac{df}{dx} \right\}(y).
\enqas \enp

\section{Examples of $A$-transform}
\subsection{Monomials.}
Appell polynomials $Q_n^\nu(y)$ are traditionally defined as
\beq
Q_n^\nu(y) = \left.\frac{d^n}{d u^n} \left(\frac{e^{uy}}{\E\left(e^{u \nu}\right)} \right)\right|_{u=0}
\enq
in other words, $\frac{e^{uy}}{\E(e^{u\nu}) }$ is the generating function for Appell polynomials
\beq
\frac{e^{uy}}{\E(e^{u\nu}) }= \sum_{n=0}^{\infty} \frac{u^n}{n!} Q_n^{\nu}(y).
\enq
\begin{prop} The ${A}$-transform with respect to $\nu$ of the monomial $y^n$ is the corresponding Appell polynomial $Q_n^\nu(y)$.
\end{prop}
\bep
Before we proceed any further, let us introduce the necessary notation.
By $\mathcal{L}^{-1}\{g\}$ we denote the inverse bilateral Laplace transform for some function $g=g(\cdot)$.
By $\delta^{(n)}(u)$ we denote the $n'$th derivative of the delta function (see \cite{GelfandShilov}, ch.I, \S 2). More precisely,
\beq \int_{-\infty}^{\infty} \delta^{(n)}(u)\phi(u) du = (-1)^n \phi^{(n)}(0) \enq

Note, that the inverse bilateral Laplace transform of $y^n$  is the $n'$th derivative of the delta function, $$\mathcal{L}^{-1}\{y^n\}(u)= (-1)^n\delta^{(n)}(u).$$
Indeed,
$$ \int_{-\infty}^\infty \mathcal{L}^{-1}\{y^n\}(u) e^{uy}du = \int_{-\infty}^\infty (-1)^n \delta^{(n)}(u) e^{uy}du = (-1)^{2n} \left.\frac{d^n}{d u^n} \left(e^{uy} \right)\right|_{u=0}= y^n.$$

Therefore,
\beqas
\mathcal{A}^\nu\{y^n\}(y) &=& \int_{-\infty}^\infty\mathcal{L}^{-1}\{y^n\}(u)\frac{e^{uy}}{\E e^{ u {\nu}}}du \\
&=& \int_{-\infty}^\infty (-1)^n \delta^{(n)}(u)\frac{e^{uy}}{\E e^{ u {\nu}}}du= \\
&=& \left.\frac{d^n}{d u^n} \left(\frac{e^{uy}}{\E\left(e^{u \nu}\right)} \right)\right|_{u=0}\\
&=& Q_n^\nu(y)
\enqas
\enp
Thus with a slight abuse of notation we write for simplicity $Q_n^\nu(y)$ instead of $Q_{y^n}^\nu(y)$.
\subsection{Polynomials (analytic functions/formal power series).}
Assume that function $g=g(\cdot)$ is a polynomial, (analytic or a formal power series in the style of umbral calculus \cite{umbral}), then  we can show, that the ${A}$-transform  of $g$ can be represented as a linear combination (power series) of Appell polynomials.
\begin{prop} \label{AppSer}
Let
\beq \label{analytic g} g(y) = \sum_{k=0}^n c_k y^k.\enq
Then ${A}$-transform of $g$ is a linear combination in Appell polynomials with the same coefficients as $g(\cdot)$, i.e.
\beqn
\mathcal{A}^{\nu} \{ g\}(y)= \sum_{k=0}^n c_k  Q_k^{\nu}(y),
\enqn
where $Q_k^{\nu}(y)$ are Appell polynomials of order $k$ generated by the random variable ${\nu}$.
\end{prop}
\bep \beqas
\mathcal{A}^\nu\{g\}(y) &=& \int_{-\infty}^\infty \mathcal{L}^{-1}\{g\}(u)\frac{e^{uy}}{\E e^{ u {\nu}}} du =
\int_{-\infty}^\infty  \left( \sum_{k=0}^n c_k  \mathcal{L}^{-1}\{y^k\}(u)\right) \frac{e^{uy}}{\E e^{ u {\nu}}} du \\
&=&  \sum_{k=0}^n c_k\int_{-\infty}^\infty  \mathcal{L}^{-1}\{y^k\}(u) \frac{e^{uy}}{\E e^{ u {\nu}}} du =
\sum_{k=0}^n c_k\int_{-\infty}^\infty (-1)^k \delta^{(k)}(u) \frac{e^{uy}}{\E e^{ u {\nu}}} du \\ &=&
\sum_{k=0}^n  c_k \frac{d^k}{ du^k}\left. \left( \frac{e^{uy}}{\E e^{ u {\nu}}} \right)\right|_{u=0}=
\sum_{k=0}^n c_k  Q_k^{\nu}(y)
\enqas

\subsection{Linear combination of exponentials.}
Let the reward function $g$ be given by linear combination of exponentials $$g(y) = \sum_{k=0}^n c_k  e^{r_k y}.$$
One can notice that  the inverse bilateral Laplace transform of $ e^{r_k y}$  is the delta function  at $r_k$, $$\mathcal{L}^{-1}\{e^{r_k y}\}(u)= \delta (u-r_k).$$
Indeed,
$$ \int_{-\infty}^\infty \mathcal{L}^{-1}\{e^{r_k y}\}(u) e^{uy}du = \int_{-\infty}^\infty  \delta(u-r_k) e^{uy}du = \left.e^{uy} \right|_{u=r_k}= e^{r_k y}.$$
\begin{prop} \label{AppSer}
Let
\beq \label{exponential_g} g(y) = \sum_{k=0}^n c_k e^{r_k y}.\enq
Then the ${A}$-transform  of $g$ is a sum of the corresponding Esscher transforms, i.e.
\beqn
\mathcal{A}^{\nu} \{ g\}(y)= \sum_{k=0}^n c_k  \frac{e^{r_k y}}{\E e^{r_k \nu}}.
\enqn
\end{prop}
\bep \beqas
\mathcal{A}^\nu\{g\}(y) &=& \int_{-\infty}^\infty \mathcal{L}^{-1}\{g\}(u)\frac{e^{uy}}{\E e^{ u {\nu}}} du =
\int_{-\infty}^\infty  \left( \sum_{k=0}^n c_k \mathcal{L}^{-1}\{e^{r_k y}\}(u)\right) \frac{e^{uy}}{\E e^{ u {\nu}}} du \\
&=&  \sum_{k=0}^n c_k\int_{-\infty}^\infty  \mathcal{L}^{-1}\{e^{r_k y}\}(u) \frac{e^{uy}}{\E e^{ u {\nu}}} du =
\sum_{k=0}^n c_k\int_{-\infty}^\infty  \delta(u-r_k) \frac{e^{uy}}{\E e^{ u {\nu}}} du \\ &=&
\sum_{k=0}^n  c_k \left. \frac{e^{uy}}{\E e^{ u {\nu}}} \right|_{u=r_k}=
\sum_{k=0}^n c_k \frac{e^{r_k y}}{\E e^{ r_k {\nu}}}.
\enqas
\enp
\subsection{Linear combinations of exponential polynomials.}
Let the reward function $g$ be given by an exponential polynomial $$g(y) = \sum_{k=0}^n c_k y^k e^{r_k y}.$$
 Note that  the inverse bilateral Laplace transform of $y^k e^{r_k y}$  is the $k'$th derivative of the delta function  at $r_k$, $$\mathcal{L}^{-1}\{y^k e^{r_k y}\}(u)= (-1)^k\delta^{(k)}(u-r_k).$$
Indeed,
$$ \int_{-\infty}^\infty \mathcal{L}^{-1}\{y^k e^{r_k y}\}(u) e^{uy}du = \int_{-\infty}^\infty (-1)^k \delta^{(k)}(u-r_k) e^{uy}du = (-1)^{2k} \left.\frac{d^k}{d u^k} \left(e^{uy} \right)\right|_{u=r_k}= y^k e^{r_k y}.$$
Denote the $k$-th derivative in $u$ of $\frac{e^{uy}}{\E e^{u\nu}}$ at $u=a$ by $Q_k^\nu(y;a)$:
\beq
Q_k^\nu(y;a):= \left.\frac{d^k}{d u^k}\left(\frac{e^{uy}}{\E e^{ u {\nu}}}\right)\right|_{u=a}
\label{Q_k}
\enq
\begin{prop} \label{ExponentialPolynomialATransform}
Let
\beq \label{exp-polynomial-g} g(y) = \sum_{k=0}^n c_k y^k e^{r_k y}.\enq
Then the ${A}$-transform  of $g$ is the following
\beqn
\mathcal{A}^{\nu}\{g\}(y)= \sum_{k=0}^n  c_k Q_k^\nu(y;r_k),
\enqn
where functions $Q=Q_k^\nu(y;r_k)$ are defined by (\ref{Q_k}).
\end{prop}
\bep \beqas
\mathcal{A}^\nu\{g\}(y)&=& \int_{-\infty}^\infty \mathcal{L}^{-1}\{g\}(u)\frac{e^{uy}}{\E e^{ u {\nu}}} du =
\int_{-\infty}^\infty  \left( \sum_{k=0}^n c_k  \mathcal{L}^{-1}\{y^k e^{r_k y}\}(u)\right) \frac{e^{uy}}{\E e^{ u {\nu}}} du \\
&=&  \sum_{k=0}^n c_k\int_{-\infty}^\infty  \mathcal{L}^{-1}\{y^k e^{r_k y}\}(u) \frac{e^{uy}}{\E e^{ u {\nu}}} du =
\sum_{k=0}^n c_k\int_{-\infty}^\infty (-1)^k \delta^{(k)}(u-r_k) \frac{e^{uy}}{\E e^{ u {\nu}}} du \\ &=&
\sum_{k=0}^n  c_k \frac{d^k}{ du^k}\left. \left( \frac{e^{uy}}{\E e^{ u {\nu}}} \right)\right|_{u=r_k}\\&=& \sum_{k=0}^n  c_k Q_k^\nu(y;r_k).
\enqas
\enp

\section{Known examples of optimal stopping problems now solved via A-transform}
\subsection{The Novikov-Shiryaev optimal stopping problem with $g(x)= (x^+)^n$.}
In \cite{NovikovShiryaevDiscr} Novikov and Shiryaev solved the optimal stopping problem (\ref{problemset}) with
$g(x) = (x^+)^n$, $n = 1, 2,\dots $ for random walks, and in [15] Kyprianou and
Surya found the solution for L\'{e}vy processes.

Here we repeat their results with our method.

Let $X$ be a L\'{e}vy process with $X_0=0$, and the reward function $g(x)= (x^+)^n$. Then we have \beqas\eta(x) &=& \argmaxproc_{0 \leq t \leq  e_q}(x+X_t)-x \\
&=& \max_{0 \leq t \leq  e_q}(x+X_t)-x \\
&=& \max_{0 \leq t \leq  e_q}(X_t)\\
&=& \eta.\enqas

The $A$-transform of $x^n$ is an Appell polynomial of order $n$
\beqas
\mathcal{A}^\eta\left\{x^n\right\}(y) &=& \int_{-\infty}^\infty \mathcal{L}^{-1}\{x^n\}(u) \frac{e^{uy}}{\E e^{u\eta}}du\\
&=& \int_{-\infty}^\infty \delta(n,u) \frac{e^{uy}}{\E e^{u\eta}}du\\
&=& Q_n^\eta(y)
\enqas
In such a way we repeat the results of Novikov and Shiryaev, and if $\left(Q_n^\eta(x)\right)^+$ and $(x^n)^+$ are co-monotone, than the optimal stopping boundary is the positive root of the Appell polynomial $Q_n^\eta(x)$.

\subsection{The Novikov-Shiryaev optimal stopping problem with $g(x)= (x^+)^\nu$.}
In \cite{NovikovShiryaevCont} Novikov and Shiryaev solved the optimal stopping problem (\ref{problemset}) with
$g(x) = (x^+)^\nu$ when the underlying process is  a L\'{e}vy process.

Here we repeat their results with our method.

Let $X$ be a L\'{e}vy process with $X_0=0$, and the reward function $g(x)= (x^+)^\nu$. Exactly as in the previous example we obtain \beqas\eta(x) &=& \argmaxproc_{0 \leq t \leq  e_q}(x+X_t)-x \\
&=& \max_{0 \leq t \leq  e_q}(x+X_t)-x \\
&=& \max_{0 \leq t \leq  e_q}(X_t)\\
&=& \eta.\enqas

The inverse bilateral Laplace transform of $x^\nu$ with $\nu<0$ is
\beqas
\mathcal{L}^{-1}\{x^\nu\}(u)= \left\{
\begin{array}{c l} \displaystyle
     -\frac{(-u)^{-\nu-1}}{\Gamma(-\nu)}, & \text{if } u < 0\\
     0,              & \text{if } u \geq 0.
\end{array}\right.
\enqas
where $\Gamma$ is a gamma function.
Indeed,
\beqas
\int_{-\infty}^\infty \mathcal{L}^{-1}\{x^\nu\}(u) e^{uy}du &=& \int_{-\infty}^0 -\frac{(-u)^{-\nu-1}}{\Gamma(-\nu)} e^{uy}du\\
&=& \int_0^\infty \frac{u^{-\nu-1}}{\Gamma(-\nu)} e^{-uy}du\\
&=& y^\nu
\enqas
Thus, for $\nu <0$
the $A$-transform with respect to $\eta$ of $x^\nu$ is given by
\beqas
\mathcal{A}^\eta\{x^\nu\}(y) &=& \int_{-\infty}^\infty \mathcal{L}^{-1}\{x^\nu\}(u) \frac{e^{uy}}{\E e^{u\eta}}du\\
&=& \int_{-\infty}^0 -\frac{(-u)^{-\nu-1}}{\Gamma(-\nu)} \frac{e^{uy}}{\E e^{u\eta}}du\\
&=&\int_{0}^\infty \frac{u^{-\nu-1}}{\Gamma(-\nu)} \frac{e^{-uy}}{\E e^{-u\eta}}du.
\enqas
This coincides  with the results obtained by Novikov and Shiryaev for $\nu<0$.
 \newline
  Therefore, we repeat the results of \cite{NovikovShiryaevCont} for $\nu<0$, and state that if $\left(\mathcal{A}^\eta\{x^\nu\}(x)\right)^+$ and $(x^\nu)^+$ are co-monotone, then the optimal stopping boundary is the positive root of the function obtained as the $A$ -transform with respect to $\eta$ of $x^\nu$, i.e. $\mathcal{A}^\eta\{x^\nu\}(x)$.
\section{New example. Two-sided problem}
Consider the optimal stopping problem (\ref{problemset})  with the reward function
\beq
g(x) = e^{a x} + e^{-b x}-2,
\enq
 with $\sqrt{2 q} >a>b>0$ some constants, and $q$ is an interest (killing) rate in our stopping problem. We assume that the underlying process is a Brownian motion $B_t$  with $B_0=0$.

 \begin{figure}[ht]
\noindent
\begin{center}
\fbox{\scalebox{.40}{\includegraphics{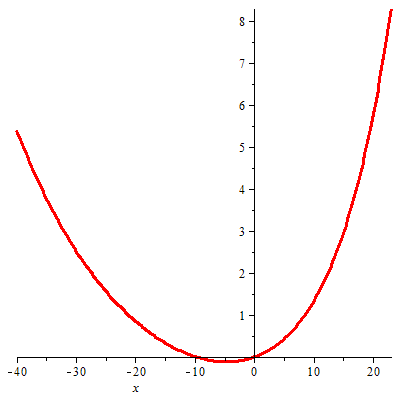}}}
\end{center}
 \caption{\footnotesize The reward function $g(x)= e^{ x/10} + e^{- x/20}-2.$
 }
\label{PayoffFunctionPicture}
\end{figure}
The function $g(x)$ is decreasing for $x \leq  \ln (b/a)/(a+b)$ and increasing for $x \geq \ln (b/a)/(a+b) $. \\Moreover, it is well known that
$ \displaystyle (\sup_{0\leq t \leq e_q}B_t) $ and $ \displaystyle (-\inf_{0\leq t \leq e_q}B_t )$ are equal in distribution.\\ To find $\eta(x)$ we  should compare
$$  g\left(x+\sup_{0\leq t \leq e_q}B_t\right) \mbox{        and         }  g\left(x+\inf_{0\leq t \leq e_q}B_t \right),$$ or, equivalently,  to compare $$  g\left(x+\sup_{0\leq t \leq e_q}B_t\right) \mbox{        and        } g\left(x-\sup_{0\leq t \leq e_q}B_t \right),$$ or, equivalently, to  compare $$\displaystyle  e^{(a+b)x}\:\: \:\mbox{ and }\:\:\:  \frac{\sinh\left(b\,  \sup_{0\leq t \leq e_q}B_t \right)}{ \sinh\left(a\,  \sup_{0\leq t \leq e_q}B_t \right)}.$$

Define function $f: [0, \infty) \rightarrow \mathbb{R}$  as \beq f(u) =\frac{ \sinh(b u)}{\sinh(a u)}.\enq
The function $f$ is decreasing on $[0,\infty)$ due to $a>b>0$. We write $f^{-1}$ for the inverse function of $f$, and denote by $c(x)$ the function $c(x) =f^{-1} \left( e^{(a+b)x}\right)$ .

It is easy to see that
  \beqa\eta(x) &=& \argmaxproc_{0 \leq t \leq  e_q}(x+B_t)-x \nonumber \\ \label{eta(x)twosided}
&=& \displaystyle \left\{\begin{array}{cc}
                 \displaystyle \sup_{0 \leq t \leq e_q}\!\!\! B_t & \mbox{ if } x \geq  \frac{ \ln (b/a)}{(a+b)}, \\
                 \displaystyle \sup_{0 \leq t \leq e_q}\! \!\! B_t & \mbox{ if }  x <  \frac{ \ln (b/a)}{(a+b)} \mbox{   and  } \displaystyle \sup_{0 \leq s \leq e_q}\!\!\! B_s \geq c(x),  \\
                 \displaystyle \inf_{0 \leq t \leq e_q}\!\!\! B_t &  \mbox{ if } x < \frac{ \ln (b/a)}{(a+b)} \mbox{   and  } \displaystyle \inf_{0 \leq s \leq e_q}\!\!\! B_s >
                 - c(x).
               \end{array} \right.
\enqa
 As our reward function $g$ is a  linear combination of exponential functions plus some constant, then $A$-transform with respect to $\eta(x)$ of  $g$ is given by
\beqas
\mathcal{A}^{\eta(x)}\{g\}(y)  = \frac{e^{a y}}{\E(e^{a \,\eta(x)})} + \frac{e^{-b y}}{\E(e^{-b\, \eta(x)})}-2.
\enqas
Let us calculate $\E(e^{u \eta(x)})$:
\beqa
\E e^{u\,\eta(x)} =\displaystyle \left\{\begin{array}{cc} \int_0^\infty e^{uy}p_{sup}(y) d y, & x \geq  \frac{ \ln (b/a)}{(a+b)},\\
\int_{c(x)}^\infty e^{uy } p_{sup}(y)d y + \int_{-c(x)}^0 e^{uy} p_{inf}(y)dy , & x < \frac{ \ln (b/a)}{(a+b)},
\end{array}
\right.
\enqa
where $p_{\inf}(y)= \sqrt{2 q}\, e^{y \sqrt{2 q}}$  and $p_{sup}(y)=\sqrt{2 q} \,e^{-y \sqrt{2 q}}$ are the probability density functions for $\inf_{0\leq t \leq e_q}B_t $ and $\sup_{0\leq t \leq e_q}B_t $ respectively.
\\
Subsequently, for $u < \sqrt{2 q}$ we have
\beqa
\E e^{u\,\eta(x)} =\displaystyle \left\{\begin{array}{cc} \frac{\sqrt{2 q}}{\sqrt{2 q} - u}, & x \geq  \frac{ \ln (b/a)}{(a+b)},\\
\frac{\sqrt{2 q}}{\sqrt{2 q} - u}e^{-c(x)(\sqrt{2 q}-u)} - \frac{\sqrt{2 q}}{\sqrt{2 q} + u} e^{-c(x)(\sqrt{2 q}+u)} + \frac{\sqrt{2 q}}{\sqrt{2 q} + u}, & x < \frac{ \ln (b/a)}{(a+b)},
\end{array}
\right.
\enqa

 Therefore  $\E e^{a \,\eta(x)}$ and  $\E e^{-b\,\eta(x)}$  are given by
\beqa
\E e^{a\,\eta(x)} =\displaystyle \left\{\begin{array}{cc} \frac{\sqrt{2 q}}{\sqrt{2 q} - a}, & x \geq  \frac{ \ln (b/a)}{(a+b)},\\
\frac{\sqrt{2 q}}{\sqrt{2 q} - a}e^{-c(x)(\sqrt{2 q}-a)} - \frac{\sqrt{2 q}}{\sqrt{2 q} + a} e^{-c(x)(\sqrt{2 q}+a)} + \frac{\sqrt{2 q}}{\sqrt{2 q} + a}, & x < \frac{ \ln (b/a)}{(a+b)},
\end{array}
\right.
\enqa
and
\beqa
\E e^{-b\,\eta(x)} =\displaystyle \left\{\begin{array}{cc} \frac{\sqrt{2 q}}{\sqrt{2 q} +b}, & x \geq  \frac{ \ln (b/a)}{(a+b)},\\
\frac{\sqrt{2 q}}{\sqrt{2 q} +b}e^{-c(x)(\sqrt{2 q}+b)} - \frac{\sqrt{2 q}}{\sqrt{2 q} -b} e^{-c(x)(\sqrt{2 q}-b)} + \frac{\sqrt{2 q}}{\sqrt{2 q} -b}, & x < \frac{ \ln (b/a)}{(a+b)},
\end{array}
\right.
\enqa

One can easily check that for each fixed $x$ the functions $g(y)$ and $\mathcal{A}^{\eta(x)}\{g\}(y)$ are co-monotone (in $y$) for those $y$ where $\mathcal{A}^{\eta(x)}\{g\}(y)$ is nonnegative. Consequently, to find the optimal stopping boundaries we have to find the zeros of $\mathcal{A}^{\eta(x)}\{g\}(x)$.

\begin{figure}[ht]
\noindent
\begin{center}
\fbox{\scalebox{.60}{\includegraphics{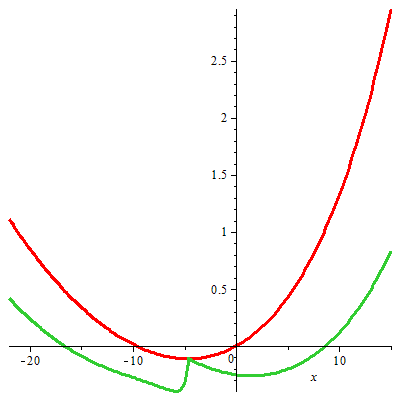}}}
\end{center}
 \caption{\footnotesize The graphs of the reward function $g=g(x)$ (red) and its  $\mathcal{A}^{\eta(x)}$-transform as a function of $x$, i.e. $A^{\eta(x)}\{g\}(x)$ (green) for $a=0.1$, $b=0.05$ and $q=0.02$.}\label{figTwosided}
\end{figure}

In other words, there are two optimal stopping boundaries $x_*$ and $x^*$, where $x^* >0$  is the zero of the equation
\beq
\frac{\sqrt{2 q} -a}{\sqrt{2 q}} e^{a \,x}+ \frac{\sqrt{2 q} +b}{\sqrt{2 q}} e^{-b \,x}- 2 =0,
\enq
and $x_*<0$  is  the zero of
\beqa \nonumber
& &\frac{e^{a x}}{\frac{\sqrt{2 q}}{\sqrt{2 q} - a}e^{-c(x)(\sqrt{2 q}-a)} - \frac{\sqrt{2 q}}{\sqrt{2 q} + a} e^{-c(x)(\sqrt{2 q}+a)} + \frac{\sqrt{2 q}}{\sqrt{2 q} + a}}
\\&+&
\frac{e^{- b x}}{\frac{\sqrt{2 q}}{\sqrt{2 q} +b}e^{-c(x)(\sqrt{2 q}+b)} - \frac{\sqrt{2 q}}{\sqrt{2 q} -b} e^{-c(x)(\sqrt{2 q}-b)} + \frac{\sqrt{2 q}}{\sqrt{2 q} -b}} -2=0,
\enqa
where $c(x) =f^{-1} \left( e^{(a+b)x}\right)$ .
\\
\\
The graph of function $\mathcal{A}^{\eta(x)}\{g\}(x)$ for $a=0.1$, $b=0.05$ and $q=0.02$ is shown in Fig.\ref{figTwosided}

\section{Conclusion and further development}

In this paper, we presented a novel approach for solving optimal stopping problems by means of applying a specially designed integral transform to the reward function. The important feature of our method that it works for non-monotone reward functions.  To construct the integral transform we need the reward function to have an inverse bilateral Laplace transform in some form.
\\
\\
The newly defined random variable $\argmaxproc$ plays the central role in the construction of the integral transform.
Calculation of $\argmaxproc$ for various Levy processes is the task to be explored.
\\
\\
Although our primary aim in this paper was to solve optimal stopping problems, it is worthwhile mentioning a by-product of our results. The integral transform we created produces a martingale if built on  a L\'{e}vy process.
\\
\\
We showed that our method
works particulary well when the reward function is a polynomial, an exponential or an exponential polynomial. This naturally leads us to explore the possibility of creating numerical methods for solving optimal stopping problems by approximating the reward functions by polynomials/exponential polynomials. The work on this topic has barely begun but looks very promising.
\\
\\
This method benefits from the straight forward generalization to multiple dimensions, which is our work in progress at the moment.
\\
\\{\bf Acknowledgements.} The author is grateful to Yuliya Mishura for fruitful discussions.

\end{document}